\documentclass[12pt]{article}
\usepackage{amsfonts}
\usepackage{}
\usepackage{wasysym,amssymb,eufrak,indentfirst,graphicx,cite,amsthm,color}
\usepackage[bookmarksnumbered, colorlinks, plainpages]{hyperref}
\newtheorem{thm}{Theorem}[section]%
\newtheorem{lemma}[thm]{Lemma}%
\newtheorem{defi}[thm]{Definition}%
\newtheorem{liz}[thm]{Example}%
\newtheorem{rem}[thm]{Remark}
\newtheorem{que}[thm]{Question}

 \def\s{\sigma}

\def\nd{\mathrel{\bigm|\kern-.7em/}}

\def\qed{\hfill $\Box$}

\begin{document}
\baselineskip 17pt

\title{On a generalisation of finite $T$-groups\thanks{Research was supported by the Fundamental Research Funds for the Central Universities (No. 2020QN20) and NNSF  of China (No. 11771409.)}}

\author{Chi Zhang\\
{\small Department of Mathematics, China University of Mining and Technology,}\\
{\small Xuzhou, 221116, P. R. China}\\
{\small E-mail: zclqq32@cumt.edu.cn}\\ \\
{ Wenbin Guo}\\
{\small  Department of Mathematics, University of Science and Technology of China,}\\
{\small Hefei, 230026, P.R. China}\\
{\small E-mail: wbguo@ustc.edu.cn}\\ \\}

\date{}
\maketitle

\begin{abstract}
Let $\s=\{\s_i |i\in I\}$ is some partition of all primes $\mathbb{P}$ and $G$ a finite group.
A subgroup $H$ of $G$ is said to be \emph{$\s$-subnormal} in $G$ if there exists a subgroup
chain $H=H_0\leq H_1\leq \cdots \leq H_n=G$ such that either $H_{i-1}$ is normal
in $H_i$ or $H_i/(H_{i-1})_{H_i}$ is a finite $\s_j$-group for some $j \in I$ for $i = 1, \ldots, n$.
We call a finite group $G$ a  \emph{$T_{\sigma}$-group} if every $\s$-subnormal subgroup is normal in $G$.

In this paper, we analyse the structure of the $T_{\sigma}$-groups and
give some characterisations of the $T_{\sigma}$-groups.

\end{abstract}

\let\thefootnoteorig\thefootnote
\renewcommand{\thefootnote}{\empty}


 \let\thefootnote\thefootnoteorig

\section{Introduction}

Throughout this paper, all groups are finite and $G$ always denotes a finite group. $\mathbb{P}$ denotes the set of all primes and
 $\pi$ denotes a set of primes.
If $n$ is an integer, then the symbol $\pi(n)$ denotes the set of all
primes dividing $n$; as usual, $\pi(G)=\pi(|G|)$, the set of all primes dividing
the order of $G$.

\subsection{$T$-groups}
A group $G$ is said to be a  \emph{$T$-group} if every subnormal subgroup of $G$ is normal in $G$.
The $T$-groups  are clearly the groups in which normality is a transitive relation.
The classical works by Gasch\"{u}tz \cite{Ga} and Robinson \cite{Rb} reveal a very detailed picture of such
groups.

Recall that $G$ is said to be a \emph{Dedekind group} if every subgroup of $G$ is normal in $G$;
It is clear that a nilpotent group $G$ is a $T$-group if and only if every subgroup
of $G$ is normal in $G$;
that is, $G$ is a Dedekind group.
More generally, Gasch\"{u}tz proved the following
result:

\begin{thm}\label{G-main}{\rm (See Gasch\"{u}tz \cite{Ga})}

Let $G$ be a group with $G^{\mathfrak{N}}$ the nilpotent residual of $G$, that is  the intersection of all
normal subgroups $N$ of $G$ with nilpotent quotient $G/N$.
Then $G$ is a soluble
$T$-group if and only if the following conditions hold:

(i) $G^{\mathfrak{N}}$ is a normal abelian Hall subgroup of $G$ with odd order;

(ii) $G/G^{\mathfrak{N}}$ is a Dedekind group;

(iii) Every subgroup of $G^{\mathfrak{N}}$ is normal in G.
\end{thm}

Recall that a group $G$ satisfies the condition
$\mathfrak{R}_p$ \cite{Rb} (where $p$ is a prime) if every subgroup of a Sylow $p$-subgroup $P$ of $G$
is normal in the normalizer of $P$.
Robinson studied the structure of finite $T$-groups using the condition $\mathfrak{R}_p$
and get the following theorem.

\begin{thm}\label{R-main}{\rm (See Robinson \cite{Rb})}

A finite group $G$ which satisfies $\mathfrak{R}_p$ for all $p$  if and only if $G$ is a soluble
$T$-group.
\end{thm}

Some other characterisations of the soluble $T$-groups have been researched.  (See \cite{Peng1, Peng2, AB2}).

\subsection{The theory of $\s$-groups}
In recent years, a new theory of $\s$-groups has been established by A. N. Skiba and W. Guo (See \cite{GS_G1, AN1, AN2, GS_China}).

In fact, following L. A. Shemetkov \cite{LA}, $\s=\{\s_i |i\in I\}$ is some partition of all primes $\mathbb{P}$, that is, $\mathbb{P}=\bigcup_{i\in I}\s_i$ and $\s_i\cap \s_j=\emptyset$ for all $i\neq j$.
$\Pi$ is always supposed to be a non-empty subset of $\sigma$ and $\Pi^{'}=\sigma\backslash \Pi$.
We write $\s(G)=\{\s_i|\s_i\cap \pi(G)\neq \emptyset\}$.

Following \cite{AN1, AN2, AN_JAA}, $G$ is said to be: \emph{$\s$-primary}
if $|\s(G)|\leq 1$; \emph{$\s$-soluble} if every chief factor of $G$ is $\s$-primary.
$G$ is called  \emph{$\sigma$-nilpotent}
 if $G = G_1 \times \cdots \times G_n$ for some $\s$-primary groups $G_1, \cdots , G_n$.
A subgroup $H$ of $G$ is said to be \emph{$\s$-subnormal} in $G$ if there exists a subgroup
chain $H=H_0\leq H_1\leq \cdots \leq H_n=G$ such that either $H_{i-1}$ is normal
in $H_i$ or $H_i/(H_{i-1})_{H_i}$ is $\s$-primary for all $i=1,\cdots,n$;
An integer $n$ is said to be a \emph{$\Pi$-number} if $\pi(n)\subseteq\bigcup_{\sigma_{i}\in \Pi}\sigma_{i}$;
a subgroup $H$ of $G$ is called a \emph{$\Pi$-subgroup} of $G$ if $|H|$ is a $\Pi$-number;
a subgroup $H$ of $G$ is called a \emph{Hall $\Pi$-subgroup} of $G$ if $H$ is a $\Pi$-subgroup of $G$ and $|G:H|$ is a $\Pi^{'}$-number;
a subgroup $H$ of $G$ is called a \emph{Hall $\s$-subgroup} of $G$ if $H$ is a Hall $\Pi$-subgroup of $G$ for some $\Pi \subseteq \s$.
We use $\mathfrak{N}_{\s}$ to denote the class of $\sigma$-nilpotent groups.

\begin{rem}
In the case when $\sigma=\{\{2\},\{3\},\cdots\}$ (we use here the notation in
\cite{AN2}),

(1) $\s$-soluble groups and $\sigma$-nilpotent groups are soluble groups and nilpotent groups respectively.

(2) $\s$-subnormal is subnormal.

(3) A  Hall $\s$-subgroup of $G$ is a Hall $\pi$-subgroup for some $\pi \subseteq \mathbb{P}$.

(4) A  Hall $\s_i$-subgroup of $G$ is a Sylow subgroup of $G$.

\end{rem}

This new theory of $\s$-groups is
actually the development and popularization of the famous Sylow theorem, the Hall theorem
of the soluble groups and the Chunihin theorem of $\pi$-soluble groups.
A series of studies have been caused (See, for example, \cite{AN1, AN2, guo-zhang-skiba, zhang-skiba-u, zss, Hu1, GS_G1, B, GS_G2, 14, AN_JAA, GS_A, GS_China, AN_C, zhang-wu-guo}).

\subsection{The $T_{\s}$-groups and the main results}

Combined with the above two contents, we naturally reposed the following problem:

\begin{que}
Based on this new theory of $\s$-groups, could we establish the theory of generalised $T$-groups $?$
\end{que}

In this paper, we will solve this question. We first give the following definition:

\begin{defi}
We called a group $G$ a $T_{\s}$-group if every $\s$-subnormal of $G$ is normal in $G$.
\end{defi}

Bearing in mind  the results in \cite{Ga, Rb},  it seems to
be natural to ask:

\begin{que}\label{Q-main}
What is the structure of the $T_{\s}$-groups $?$
\end{que}

We will give a complete answer to this question in the case when
$G$ is $\s$-soluble.
It is clear that every subnormal subgroup is $\s$-subnormal in $G$
and so every $T_{\s}$-groups is a  $T$-groups.
However, the following example shows that the converse is not true.

\begin{liz}
Let $A = C_3 \rtimes C_2$ be a non-abelain subgroup of order 6 and let $G = A \times C_5$.
Let  $\sigma=\{\sigma_1,\sigma_2, \s_3 \}$, where $\sigma_1=\{2,3\}$, $\s_2= \{ 5 \}$ and $\sigma_3=\{2,3,5\}'$.
Then $G$ is a $T$-group but is not a $T_{\s}$-group. In fact, obviously, $G$ is a $T$-group.
However $G$ is  not a  $T_{\s}$-group since $C_2$ is a $\s$-subnormal subgroup of $G$
but is not normal in $G$.
\end{liz}

In order to better describe the $T_{\s}$-groups, I give the following definition:

\begin{defi}
We called a group $G$ satisfies the condition $\mathfrak{R_{\s_i}}$ if  every subgroup $K$ of every Hall $\pi$-subgroup $H$ of $G$ (for $\pi \subseteq \s_i$) is normal in the normalizer $N_G(H)$ of $H$.
\end{defi}

\begin{rem}
In the case when $\sigma=\{\{2\},\{3\},\cdots\}$, the condition $\mathfrak{R_{\s_i}}$ is just the condition $\mathfrak{R}_p$.
\end{rem}

The following theorem gives a answer to Question \ref{Q-main}.

\begin{thm}\label{main}

Let $G$ be a group, $D = G^{\mathfrak{N}_{\s}}$ and $G$ is  $\s$-soluble.
Then the following statements are equivalent.

(1) $G$ is a $T_{\s}$-group;

(2) $G$ satisfies the conditions $\mathfrak{R_{\s_i}}$ for all $i$.

(3) $G$ satisfies the following conditions:

~~(i) $G = D \rtimes M$, where $D$ is an abelian Hall subgroup of $G$ of odd order, $M$ is a Dedekind group;

~~(ii) every element of $G$ induces a power automorphism on $D$; and

~~(iii) $O_{\s_i}(D)$ has a normal complement in a Hall $\s_i$-subgroup of $G$ for all $i$.
\end{thm}

In this theorem, $G^{\mathfrak{N}_{\s}}$ denotes the $\s$-nilpotent residual of $G$, that is,  the intersection of all
normal subgroups $N$ of $G$ with $\s$-nilpotent quotient $G/N$, and $O_{\s_i}(D)$ denotes the maximal normal $\s_i$-subgroup.

\begin{rem}
In the case when $\sigma=\{\{2\},\{3\},\cdots\}$, The Theorems \ref{G-main} and \ref{R-main} are  corollaries of our Theorem \ref{main}.
\end{rem}

\section{Preliminaries}

\begin{lemma}\label{nilp}{\rm\cite[Corollary 2.4 and Lemma 2.5]{AN1}}
The class $\mathfrak{N}_{\s}$ of all $\sigma$-nilpotent groups is closed under taking products of normal subgroups, homomorphic images and subgroups.
Moreover, if $E$ is a normal subgroup of $G$ and $E/E\cap \Phi(G)$ is $\sigma$-nilpotent, then $E$ is $\sigma$-nilpotent.
\end{lemma}

\begin{lemma}\label{subgroup}{\rm\cite[Lemma 2.6(6)]{AN1}}
every subgroup of a $\s$-nilpotent group $G$ is $\s$-subnormal in $G$.
\end{lemma}

The following lemma directly follows from Lemma \ref{nilp} and \cite[Lemma 1.2]{LA} (see also \cite[Chap. 1, Lemma 1.1]{GuoII}).

\begin{lemma}\label{re}
If $N$ is a normal subgroup of $G$, then
$(G/N)^{\mathfrak{N}_{\s}} = G^{\mathfrak{N}_{\s}}N/N$.
\end{lemma}

\begin{lemma}\label{subnormal}{\rm\cite[Lemma 2.6]{AN1}}
Let $A,K$ and $N$ be subgroups of $G$. Suppose that $A$ is $\sigma$-subnormal in $G$ and $N$ is normal in $G$.
Then:

$(1)$ $A\cap K$ is $\sigma$-subnormal in $K$.

$(2)$ $AN/N$ is $\sigma$-subnormal in $G/N$.

$(3)$ If $K$ is a $\s$-subnormal subgroup of $A$, then $K$ is $\s$-subnormal in $G$.

$(4)$ If $A$ is a $\s$-Hall subgroup of $G$, then $A$ is normal in $G$.

$(5)$ If $H \neq 1$ is a Hall $\s_i$-subgroup of $G$ and
$A$ is not a $\s_i^{'}$-group, then $A \cap H \neq 1$ and
$A \cap H$ is a Hall $\s_i$subgroup of $A$.
\end{lemma}

\begin{lemma}\label{hall}{\rm(P. Hall\cite{hall})}
Let $G$ be a soluble group and $\pi$ a set of primes.
Then:

$(1)$ Hall $\pi$-subgroups of $G$ exists,

$(2)$ they form a conjugacy class of $G$, and

$(3)$ each $\pi$-subgroup of $G$ is contained in a Hall $\pi$-sybgroup of $G$.
\end{lemma}

The following lemma is clear.

\begin{lemma}\label{De}

(i) Every Dedekind group is nilpotent.

(ii) If $G = A \times B$, where $A$ is a Hall subgroup of $G$ and $A$ and $B$ are Dedekind groups,
then $G$ is a Dedekind group.

(iii) Every subgroup and every quotient of a Dedekind group is a Dedekind group.
\end{lemma}

\section{Proof of Theorems \ref{main}}

{\bf\large $(1) \Longrightarrow (2)$:}

Suppose that $G$ is a $T_{\s}$-group and $K$ is a subgroup of a $\pi$-Hall subgroup $H$ of $G$, where $\pi \subseteq \s_i$ for some $i$.
Since $H$ is a $\s_i$-group, $K$ is $\s$-subnormal in $H$ by Lemma \ref{subgroup}.
Note that $H$ is normal in $N_G(H)$. It implies that $K$ is $\s$-subnormal in $N_G(H)$ by Lemma \ref{subnormal}(3).
Hence $K$ is normal in $N_G(H)$ by the hypothesis. Consequently, $G$ satisfies the condition $\mathfrak{R}_{\s_i}$.

{\bf\large $(2) \Longrightarrow (3)$:}

Assume that this is false and let $G$ be a counterexample of minimal order.
We proceed via the following steps.

$(1)$ {\sl Every Hall $\s_i$-subgroup of $G$ is a Dedekind group for all $i$. Hence $D\neq 1$ and $G$ is soluble.}

Let $H$ be a Hall $\s_i$-subgroup and $K$ be a subgroup of $H$.
Then $K$ is normal in $N_G(H)$ by the hypothesis and so $K$ is normal in $H$.
Hence $H$ is a Dedekind group.
This implies that $D\neq 1$.
We now show that $G$ is soluble.
In fact, since $G$ is $\sigma$-soluble, every chief factor $S/K$ of $G$ is $\sigma$-primary, that is, $S/K$ is a $\sigma_{i}$-group for some $i$.
But as every every Hall $\s_i$-subgroup of $G$ is a Dedekind group, every Hall $\s_i$-subgroup  is  nilpotent.
Hence $S/K$ is a elementary abelian group.
It follows that $G$ is soluble.

$(2)$ {\sl Let $R$ be a non-identity minimal normal subgroup of $G$. Then
the hypothesis holds for $G/R$.
Hence $G/R$ satisfies statement $(3)$ of the Theorem.}

Let $H/R$ be a Hall $\pi$-Hall subgroup of $G/R$ where $\pi \subseteq \s_i$ and $K/R$ is a subgroup of $H/R$.
Note that $R$ is a $p$-group since $G$ is soluble by Claim (1).
Assume that $p$ belongs to $\pi$, then $H$ is a Hall $\pi$-subgroup of $G$.
Hence $K$ is normal in $N_G(H)$ by hypothesis.
Then $K/R$ is normal in $N_G(H)/R = N_{G/R}(H/R)$.
If  $p$ does not belong to  $\pi$, then there are a Hall $\pi$-subgroup $V$ of $K$ and a Hall $\pi$-subgroup $W$ of $H$
such that $V \leqslant W$ by Lemma \ref{hall}.
It is clear that $W$ is also a Hall $\pi$-subgroup of $G$ since $H/R$ be a Hall $\pi$-Hall subgroup of $G/R$.
Hence  $V$ is normal in $N_G(W)$ by hypothesis
and so $K/R = VR/R$ is normal in $N_G(W)R/R = N_{G/R}(WR/R) =N_{G/R}(H/R)$

$(3)$ {\sl The hypothesis holds for every proper Hall subgroup $M$ of $G$ and $M^{\mathfrak{N}_{\s}} \leq D$.}

Let $M_i$ be a Hall $\s_i$-subgroup of $M$  and $K$ is a subgroup of $M_i$ for all $i$.
Then $M_i$ is a Hall $\pi$-subgroup of $G$ where $\pi \subseteq \s_i$
since $M$ is a Hall subgroup.
Hence $K$ is normal in $N_G(M_i)$, and so $K$ is normal in $N_M(M_i)$.
Therefore $M$ satisfies the condition $\mathfrak{R}_{\s_i}$ for all $i$.
This shows that the hypothesis for $M$.
Moreover, since $G/D \in \mathfrak{N}_{\s}$ and $\mathfrak{N}_{\s}$ is subgroup closed by
Lemma \ref{nilp},
$$M/M \cap D \cong MD/D \in \mathfrak{N}_{\s}.$$
Hence $M^{\mathfrak{N}_{\s}} \leq M \cap D \leq D$.

$(4)$ {\sl $D$ is nilpotent}.

Assume that this is false and let $R$ be a minimal normal subgroup of $G$. Then:

$(a)$ {\sl $R = C_G(R) = O_p(G) = F(G) \leq D$ for some $p \in \s_i$. Hence $R$ is an unique minimal normal subgroup of $G$ and $R$ is not cyclic.}

First note that $RD/R = (G/R)^{\mathfrak{N}_{\s}}$ is abelian by Lemma \ref{re} and Claim $(2)$.
Therefore $R \leq D$, and so $R$ is the unique minimal normal subgroup of $G$ and $R \nleq \Phi(G)$ by Lemma \ref{nilp}.
Since $G$ is soluble, $R$ is an abelian subgroup.
It follows that $R = C_G(R) = O_p(G) = F(G)$ for some $p \in \s_i$ by \cite[Chap. A, 13.8(b)]{Doerk}.
If $|R| = p$,
then $G/R = C_G(R)$ is cyclic and so $G$ is supersoluble.
But then $D = G^{\mathfrak{N}_{\s}} \leq G' \leq F(G)$
and so $D$ is nilpotent, a contradiction. Thus $R$ is not cyclic.

$(b)$ {\sl Every Hall $\s_i$-subgroup is a Sylow $p$-subgroup, where $p \in \s_i$.}

Let $H$ be a Hall  $\s_i$-subgroup and $p \in \s_i$.
Then $H$ is nilpotent and $R \leq H$ by Claim $(1)$ and Lemma \ref{De}.
Hence  $H$ is a Sylow  $p$-subgroup by Claim $(a)$.

$(c)$ {\sl $|\pi(G)| = 2$.}

Let $H_i$ be a Hall $\s_i$-subgroup of $G$, where for $p \in \s_i$.
By Claim $(b)$, $H_i$ is a Sylow $p$-subgroup denoted by $P$.
If $|\pi(G)| = 1$, then $G$ is nilpotent, a contradiction.
Assume that $|\pi(G)| \geq 3$.
Then there exist two different primes belonging to $p'$, denoted by $q$ and $t$.
Since $G$ is soluble by Claim $(1)$,
there are a Hall $t'$-subgroup $M_1$ of $G$ and a Hall $q'$-subgroup $M_2$ of $G$ by Lemma \ref{hall}(1).
Let $V_1 = M_1^{\mathfrak{N}_{\s}}$  and  $V_2 = M_2^{\mathfrak{N}_{\s}}$.
Suppose that $V_1 = 1$ or $V_2 = 1$.
Assume without of generality that $V_1 = 1$.
Then $M_1$ is $\s$-nilpotent.
Let $Q$ be a Sylow $q$-subgroup of  $M_1$.
Then $Q \leq C_G(R) = R$ since $M_1$ is $\s$-nilpotent.
This contradiction shows that $V_1 \neq 1$ and $V_2 \neq 1$.
Since $M_1$ and $M_2$ are Hall subgroups of $G$, $M_1$ and $M_2$ satisfy the conditions $\mathfrak{R}_{\s_i}$ for all $i$ by Claim $(3)$.
The choice of $G$ implies that $V_1$ and $V_2$ are abelian Hall subgroups of $G$.
Then $R \leq V_1 \cap V_2$.
In fact, if  $R \nleq V_1$, then $R \cap V_1 = 1$.
It follows that $V_1 \leq C_{M_1}(R) = R$, and so $R = V_1$, a  contradiction.
Hence $R \leq V_1$.
Similarly, we have  $R \leq V_2$.
Note that $R$ is not cyclic by Claim $(a)$.
Let $L < R$ and $|L| = p$.
By Claim $(3)$ and the choice of $G$, every element of $M_i$ $(i = 1, 2)$ induces a power automorphism on $V_i$.
Hence $L$ are normal in  $M_1$ and $M_2$.
It follows that $L$ is normal in $\langle M_1, M_2 \rangle = G$, a contradiction.
Hence we have Claim $(c)$.

(d) {\sl The final contradiction for Claim $(4)$.}

By Claim $(c)$, we may assume that $G = PQ$ where $P$  is a Sylow $p$-subgroup of $G$ and $Q$ is a $q$-subgroup of $G$.
Since every Dedekind group of odd order is abelian by \cite[Theorem 5.3.7]{R}, we have that either $P$ is abelian or $Q$ is abelian.
If $P$ is abelian, then $RP =P$ is normal in $G$  by Claim $(a)$ and Theorem 3.2.28 in \cite{AB2}.
Hence  $D \leq P$ is nilpotent, a contradiction.
If $Q$ is abelian, then $RQ$ is normal in $G$  by Claim $(a)$ and Theorem 3.2.28 in \cite{AB2}.
Hence by Frattini augument, $G = RQN_G(Q) =  RN_G(Q)$.
Let $N_p$ is a  Sylow $p$-subgroup of $N_G(Q)$.
If $N_p = 1$, then $R$ is a normal Sylow $p$-subgroup of $G$.
Therefore $D \leq R$ is nilpotent, a contradiction.
Assume that $N_p \neq 1$.
Since  $RN_p$ is a Dedekind subgroup of $G$ by Claim (1),
$R \leq RN_p \leq N_G(N_p).$
But since $R$ is the unique minimal normal subgroup of $G$ by Claim $(a)$,
we have
$$R \leq {N_p}^G = {N_p}^{RN_G(N_p)} = {N_p}^{N_G(Q)} \leq N_G(Q).$$
It follows that $G = N_G(Q)$.
Then $Q$ is normal in $G$.
Therefore $D \leq Q$ is nilpotent.
This contradiction shows that Claim $(4)$ holds.

$(5)$ {\sl $D$ is a Hall subgroup of $G$.}

Assume that this is false. Let $P$ be a Sylow $p$-subgroup of $D$ such that $1<P<G_p$ for some prime $p$ and some Sylow
$p$-subgroup $G_p$ of $G$. Then $p||G:D|$.
We can assume without loss of generality that $G_p \leq H_1$, where $H_1$ is a Hall $\s_1$-subgroup of $G$.

$(a')$ {$D=P$ is a minimal normal subgroup of $G$.}

Let $R$ be a minimal subgroup of $G$ contained in $D$.
Then by Claim $(4)$, $R$ is a $q$-group for some prime $q$.
Moreover, $D/R = (G/R)^{\mathfrak{N}_{\s}}$ is a Hall
subgroup of $G/R$ by Claim $(2)$ and Lemma \ref{re}.
Suppose that $PR/R \neq 1$.
Then $PR/R \in Syl_p(G/R)$.
If $q \neq p$, then $P \in Syl_p(G)$.
This contradicts the fact that $P < G_p$.
Hence $q = p$,
so $R \leq  P$ and $P/R \in Syl_p(G/R)$.
We again get that $P \in Syl_p(G)$.
This contradiction shows that $PR/R = 1$,
which implies that $R = P$ is the unique
minimal normal subgroup of $G$ contained in $D$.
Since $D$ is nilpotent by Claim $(4)$,
a $p'$-complement $E$ of $D$ is characteristic in $D$
and so it is normal in $G$.
Hence $E = 1$, which implies that $R = D = P$.

$(b')$ {\sl $D\nleq \Phi(G)$. Hence for some maximal subgroup $M$ of $G$ we have $G=D\rtimes M$.}

Note that $G/D = G/G^{\mathfrak{N}_{\s}}$ is $\s$-nilpotent.
If $D \leq \Phi(G)$, then  by lemma \ref{nilp}, $G \in \mathfrak{N}_{\s}$ and so $D = 1$, a contradiction.

$(c')$ {\sl Let $R$ be a minimal normal subgroup of $G$. If $D\neq R$, then $G_p=D\times (G_p\cap R)$. Hence $O_{p'}(G)=1$.}

By Claims $(2)$ and $(a')$, we have that $DR/R$ is a Sylow $p$-subgroup of $G/R$.
It follows that $DR/R=G_pR/R$. Hence $G_p=D\times (G_p\cap R)$. Thus $O_{p'}(G)=1$ since $G$ is soluble and $D<G_p$ by Claim $(a')$.

$(d')$ {\sl Let $V=C_G(D)\cap M$. Then $V\unlhd G$ and $C_G(D)=D\times V \leq H_1$.}

In view of Claims $(a')$ and $(b')$, we have that $C_G(D)=D\times V$ and $V$
is a normal subgroup of $G$.
Moreover, $V\cong VD/D$ is $\s$-nilpotent by Lemma \ref{nilp}.
Let $W$ be a $\s_1$-complement of $V$ .
Then $W$ is characteristic in $V$ and so it is normal in $G$.
Then $W = 1$ by Claim $(c')$.
Hence we have Claim $(d')$.

$(e')$ {\sl $H_1 = G_p$ is a Sylow $p$-subgroup of $G$.}

Since $G/D$ is $\s$-nilpotent and $D \leq H_1$ by Claim $(a')$, $H_1$ is normal in $G$.
A $p'$-complement $E$ of $H_1$ is characteristic in $H_1$ since $H_1$ is nilpotent by Claim $(1)$.
Hence $E =1$ by Claim $(c')$.
It follows that $H_1 = G_p$ is a Sylow $p$-subgroup of $G$.

$(f')$ {\sl $|\pi(G)| = 2$.}

If $|\pi(G)| = 1$, then $G$ is nilpotent, a contradiction.
Assume that $|\pi(G)| \geq 3$.
Then there exist one more primes belonging to $p' = \s_1^{'}$ and let $q \in p'$.
Since $G$ is soluble, $G$ has a Hall $\{p, q\}$-subgroup $H$ of $G$.
Let  $L = H^{\mathfrak{N}_{\s}}$ and let $H = G_pQ$ where $Q$ is a $q$-subgroup of $G$.
Note that $H < G$.
If $L = 1$, then $H = P \times Q$ by Claim $(e')$.
Consequently, $Q \leq C_G(D) \leq H_1 = G_p$ by Claim $(d')$, a contradiction.
Hence $L \neq 1$.
By Claim $(3)$ and Claim $(a')$, $L \leq D = P$ and $L$ is a Hall subgroup of $H$ by the choice of $G$.
Note that $L$ is a Hall subgroup of $G$ since $H$ is a Hall subgroup of $G$.
Therefore $L = D$ is a Hall subgroup of $G$.
The contradiction shows that Claim $(f')$ holds.

$(g')$ {\sl The final contradiction for Claim $(5)$.}

Let $\pi(G) = \{p, q\}$.
Then $G$ satisfies the conditions $\mathfrak{R}_p$ and $\mathfrak{R}_q$ by Claim $(e')$ and Claim $(3)$.
Hence $G$ is a $T$-group by Theorem \ref{R-main}.
By Claims $(e')$ and $(f')$, it is clear that $G^{\mathfrak{N}} = D$.
Hence by Theorem \ref{G-main}, $D$ is a Hall subgroup of $G$.
The contradiction completes the proof of Claim $(5)$.

$(6)$ {\sl $G = D \times M$ where $M$ is a Dedekind group.}

Since $D$ is a normal subgroup of $G$, by Schur-Zassenhaus Theorem,
$G = D \times M$ and $M$ is a Hall subgroup of $G$.
But since $D = G^{\mathfrak{N}_{\s}}$, $M$ is $\s$-nilpotent.
Then by Claim $(1)$ and Lemma \ref{De}$(ii)$, we have that $M$ is a Dedekind group.

$(7)$ {\sl Let $H_i$ be a Hall $\s_i$-subgroup of $G$ for each $\s_i \in \s(D)$. Then $H_i = O_{\s_i}(D) \times S$ for some subgroup $S$ of $H_i$.}

By Claim $(1)$, $H_i$ is nilpotent.
By Claims $(4)$ and $(5)$, $D$ is a nilpotent Hall subgroup of $G$.
Hence we have Claim $(7)$

$(8)$ {\sl Every subgroup $H$ of $D$ is normal in $G$. Hence every element of $G$ induces a power automorphism on $D$.}

Since $D$ is nilpotent by Claim $(4)$,
it is enough to consider the case when $H \leq O_{\s_i}(D) = H_i \cap D$ for some $\s_i \in \s(D)$.
By condition $(2)$, $H$ is normal in $N_G(O_{\s_i}(D))$.
But clearly $N_G(O_{\s_i}(D)) = G$.
Therefore $H$ is normal in $G$.

$(9)$ {\sl $|D|$ is odd.}

Suppose that $2$ divides $|D|$. Then by Claims $(4)$ and $(7)$, $G$ has a chief factor $D/K$ with
$|D/K|= 2$.
This implies that $D/K\leq Z(G/K)$.
Since $D$ is a  normal Hall subgroup of $G$ by Claim $(5)$, it has a complement $M$ in $G$.
Hence $G/K = D/K \times MK/K$, where $MK/K \cong M \cong G/D$ is $\s$-nilpotent.
Therefore $G/K$ is $\s$-nilpotent by Lemma \ref{nilp} and Claim $(4)$.
But then $D \leq K <D$, a contradiction. Hence we have (9).

$(10)$ {\sl $D$ is abelian.}

By Claim $(8)$,  $D$ is a Dedekind group.
But $D$ is odd order by Claim $(9)$. Hence $D$ is abelian by \cite[Theorem 5.3.7]{R}.

$(11)$ {\sl Final contradiction.}

Claims $(5)$, $(6)$, $(7)$, $(8)$, $(9)$ and $(10)$ show that the conclusion (3) holds for $G$.
This final contradiction completes the proof of $(2) \Longrightarrow (3)$.

{\bf\large $(3) \Longrightarrow (1)$:}

Suppose that $G$ satisfies the conditions $(i)$, $(ii)$ and $(iii)$ of $(3)$.
Then $G$ is soluble.
Now we need to prove that every $\s$-subnormal subgroup $H$ of $G$ is normal in $G$.
Suppose that this is false, that is, some $\s$-subnormal subgroup $H$ of $G$ is not normal
in $G$.
Let $G$ be a counterexample with $|G| + |H|$ minimal.
Then by the condition $(i)$ and Lemma \ref{De}$(i)$,
we see that $D \neq 1$.
We now proceed the proof via the following steps.

$(I)$ {\sl The hypothesis holds for every quotient $G/N$ of $G$, where $N$ is a proper normal subgroup of $G$.}

By the condition $(i)$, we have that $G/N = (DN/N) \rtimes (MN/N)$, where $DN/N \cong D/D \cap N$ is an abelian
Hall subgroup of $G/N$ of odd order and $MN/N \cong M/M \cap N$ is a Dedekind-group by
Lemma \ref{De}$(iii)$. Hence condition $(i)$ holds for $G/N$. Suppose that $V/N$ is any subgroup
of $DN/N$, then $V = N(D \cap V)$.
Since $D \cap V$ is normal in $G$ by condition $(ii)$, $V/N$ is normal in $G/N$.
Hence the condition $(ii)$ holds for $G/N$.
Since $D$ is nilpotent, clearly $O_{\s_i}(D)N/N = O_{\s_i}(DN/N)$.
Condition $(iii)$ implies that $O_{\s_i}(D)$ has a normal complement $S$ in a Hall $\s_i$-subgroup $E$ of $G$ for every $i$.
Then $EN/N$ is a Hall $\s_i$-subgroup of $G/N$ and $SN/N$ is
normal in $EN/N$.
Hence $$(SN/N)(O_{\s_i}(DN/N)) = (SN/N)(O_{\s_i}(D)N/N) = EN/N$$
and
$$(SN/N) \cap O_{\s_i}(DN/N) = (SN/N) \cap (O_{\s_i}(D)N/N) = N(S \cap O_{\s_i}(D)N)/N$$
$$ = N(S \cap O_{\s_i}(D))(S \cap N)/N = N/N.$$
Hence condition $(iii)$ also holds on $G/N$.

$(II)$  {\sl $H_G = 1$.}

Assume $H_G \neq 1$.
The hypothesis holds for $G/H_G$ by Claim $(I)$.
On the other hand,
$H/H_G$ is $\s$-subnormal in $G/H_G$ by Lemma \ref{subnormal}$(2)$, so $H/H_G$ is normal in $G/H_G$ by the choice of $G$.
But then $H$ is normal in $G$, a contradiction.
Hence we have Claim $(II)$.

$(III)$ {\sl $H$ is a $\s_i$-group for some $i$ and $H \in M^x$ for all $x \in G$.}

Claim $(II)$ and the condition $(ii)$ imply that $H \cap D = 1$.
Since $H \cong HD/D \leq G/D$, $H$ is $\s$-nilpotent by Lemma \ref{nilp}.
Hence $H = A_1 \times \cdots \times A_n$ for some $\s$-primary groups $A_1, \cdots ,A_n$.
Then $H = A_i$ is a $\s_i$-group for some $i$ since otherwise $H$ is normal in $G$ by the choice of $(G, H)$.
Note that $G = D \times M$ by the condition $(i)$.
Let $M_i$ be the Hall $\s_i$-subgroup of $M$ and $E$ be a Hall $\s_i$-subgroup of $G$ containing $M_i$.
Lemma \ref{subnormal}$(5)$ implies that $H \leq E^x$ for all $x \in G$.
If $E \cap D = 1$, then $M_i$ is a Hall $\s_i$-subgroup of $G$, and so $H \leq M^x$ for all $x \in G$.
Now suppose that $E \cap D \neq 1$.
Then $H \leq E^x = O_{\s_i}(D) \times M_i^{x}$ by condition $(iii)$.
But since $H \cap D = 1$, we have also that $H \leq M_i^x \leq M^x$ for all $x \in G$.

$(IV)$ {\sl The Hall $\s_j$-subgroups of $G$ are Dedekind-groups for all $j$.}

Let $A$ be a Hall $\s_j$-subgroup of $G$.
If $A \cap D = 1$,
then $A \cong AD/D \leq G/D$, where $G/D$ is a Dedekind group by the condition $(i)$.
Hence $A$ is a Dedekind group by Lemma \ref{De}$(iii)$.
Now assume that $A \cap D \neq 1$.
Then $A = (A \cap D) \times S$ by condition $(iii)$, where $A \cap D = O_{\s_j}(D)$ and $S$ is a normal complement of $A \cap D$ in  $A$.
Then $A$ is a Dedekind group by Lemma \ref{De}$(ii)$ because $A \cap D$ and $S \cong DS/D \leq G/D$ are Dedekind groups.

$(V)$ {\sl $D$ is also a $\s_i$-group.}

Assume that this is false.
Note that $D$ is an abelian group.
Assume without of generality that $O_{\s_j}(D) \neq 1$, where $j \neq i$.
Then by Claim $(1)$,
$HO_{\s_j}(D)/O_{\s_j}(D)$ is normal in  $G/O_{\s_j}(D)$, and so  $HO_{\s_j}(D)$ is normal in $G$.
But by Lemma \ref{subnormal}(1), $H$ is $\s$-subnormal in $HO_{\s_j}(D)$.
Hence by Lemma \ref{subnormal}(4), $H$ is normal in $HO_{\s_j}(G)$.
Then $H$ is characteristic in $HO_{\s_j}(G)$.
It follows that $H$ is normal in $G$,
a contradiction.
Hence we have Claim $(V)$.

$(VI)$ {\sl Final contradiction.}

Since $D$ is a $\s_i$-group by Claim $(V)$ and $G=D\rtimes M$ by the condition (i), we have that $G = H_iM$, where $H_i$ is a Hall $\s_i$-group of $G$.
Since $G$ is soluble and $H$ is a $\s_i$-group by Claim $(III)$, we can assume without of generality that $H \leq H_i$ by Lemma \ref{hall}$(3)$.
Then $H_i \leq N_G(H)$ by Claim $(IV)$.
On the other hand,
since $M$ is a Dedekind group by the hypothesis, we have $M \leq N_G(H)$ by Claim $(III)$.
Hence $G = H_iM \leq N_G(H)$.
This shows that $H$ is normal in $G$.
This contradiction completes the proof for $(3) \Longrightarrow (1)$.

In summary, the Theorem \ref{main} is proved.

\qed


\smallskip




\begin{thebibliography}{00}

\bibitem{14}  Kh. A. Al-Sharo, A. N. Skiba, On finite groups
with $\sigma$-subnormal Schmidt subgroups, Comm. Algebra {\bf 45} (2017), 4158--4165.

\bibitem{AB2}
A. Ballester-Bolinches, R. Esteban-Romero, M. Asaad, \textit{Products of Finite Groups}, Walter de Gruyter, Berlin-New York, 2010.

\bibitem{B} J.C. Beidleman,  A. N. Skiba, On $\tau_{\sigma}$-quasinormal subgroups of finite groups,  J. Group Theory  {\bf
20}(5)  (2017), 955-964.


\bibitem{Doerk}
K. Doerk, T. Hawkes, \textit{Finite Soluble Groups}, Walter de Gruyter, Berlin, 1992.


\bibitem{Ga}
W. Gasch\"{u}tz, Gruppen, in denen das Normalteilersein transitivist,  J. Reine Angew. Math. \textbf{198} (1957), 87-92.

\bibitem{GuoII}
W. Guo, \textit{Structure Theory for Canonical Classes of Finite Groups}, Springer, Heidelberg-New York-Dordrecht-London, 2015.


\bibitem{GS_G1}  W. Guo, A. N. Skiba, Finite groups with permutable  complete
  Wielandt  sets of subgroups, J. Group Theory {\bf 18} (2015),  191--200.


\bibitem{GS_G2} W. Guo, A. N. Skiba, Groups with maximal subgroups of Sylow subgroups
 $\sigma$-permutably embedded, J. Group Theory {\bf 20}(1)
(2017), 169--183.


\bibitem{GS_A} W. Guo, A. N. Skiba, On the lattice of $\Pi_{\mathfrak{I}}$-subnormal
 subgroups of a  finite group, Bull. Austral. Math. Soc. {\bf 96}(2) (2017), 233-244.


\bibitem{GS_China} W. Guo,  A. N.  Skiba, Finite groups whose $n$-maximal subgroups
 are $\sigma$-subnormal, Science China Mathematics {\bf 62}(7) (2019), 1355-1372.


\bibitem{guo-zhang-skiba}
W. Guo, C. Zhang, A. N. Skiba, On $\s$-supersoluble groups and one generalization of $CLT$-groups, J. Algebra \textbf{512} (2018),
92-108.


\bibitem{hall}
P. Hall, Theorem like Sylow's,  Proc. London Math. Soc. \textbf{6}(3) (1956), 286-304.

\bibitem{Hu1} J. Huang, B. Hu, A. N. Skiba,   A generalisation of finite PT-groups, Bull. Aust. Math. Soc. {\bf 97}(3)
 (2018), 396-405.



\bibitem{Peng1}
T. A. Peng, Finte groups with pronormal subgroups, Proc. Amer. Math. Soc. \textbf{20} (1969), 232-234.

\bibitem{Peng2}
T. A. Peng, Pronormality in finite groups, J. London Math. Soc.  \textbf{3}(2) (1971), 301-306.

\bibitem{R}
D. J. S. Robinson, \textit{A Course in the Theory of Groups}, Springer-Verlag, Heidelberg-New York-Berlin, 1982.

\bibitem{Rb}
 D. J. S. Robinson,  A note on finite groups in which normality is transitive,
 Proc. Amer. Math. Soc. \textbf{19} (1968), 933-937.


\bibitem{LA}
 L. A. Shemetkov, \textit{Formation of Finite Groups}, Nauka, Main Editorial Board for Physical and Mathematical Literature, Moscow, 1978.



\bibitem{AN1}
A. N. Skiba, On $\sigma$-subnormal and $\sigma$-permutable subgroups of finite groups, J. Algebra \textbf{436} (2015),
1-16.


\bibitem{AN2}
 A. N. Skiba, Some characterizations of finite $\sigma$-soluble $P\sigma T$-groups, J. Algebra \textbf{495} (2018),
114-129.


 \bibitem{AN_JAA}   A. N. Skiba, A generalization of a Hall theorem, J. Algebra Appl.
 {\bf 15}(4) (2015), 21--36.


\bibitem{AN_C}  A. N. Skiba, On some results in the theory of finite partially soluble groups,
Comm. Math. Stat. {\bf 4} (2016), 281--309.



\bibitem{zhang-skiba-u}
C. Zhang, A. N. Skiba,  On $\sum_t^{\s}$-closed classes of finite groups, Ukrainian Math. J. \textbf{70}(12) (2018),
1707-1716.

\bibitem{zss}
C. Zhang, V. G. Safonov, A. N. Skiba, On $n$-multiply $\s$-local formations of finite groups, Comm.
Algebra \textbf{47}(3) (2019), 957-968.




\bibitem{zhang-wu-guo}
C. Zhang, Z. Wu and W. Guo, On weakly $\sigma$-permutable subgroups of finite groups, Publ. Math. Debrecen \textbf {91} (2017), 489-502.






\end{thebibliography}
\end{document}